\documentclass[12pt,psamsfonts]{article}

\usepackage{theorem}
\usepackage{amssymb}
\usepackage{amsmath}
\usepackage{amsfonts}

\newtheorem{thm}{Theorem}[section]
\newtheorem{prop}[thm]{Proposition}
\newtheorem{cor}[thm]{Corollary}
\newtheorem{lm}[thm]{Lemma}
\newtheorem{conj}[thm]{Conjecture}

{\theorembodyfont{\rm}

}

\def\ni{\noindent}

\def\l{\lambda}

\def\ra{\rightarrow}

\def\Sym{{\rm Sym}}
\def\pic{{\rm Pic}}
\def\ord{{\rm ord}}

\def\C{{\Bbb C}}

\def\G{{\Bbb G}}
\def\P{{\Bbb P}}
\def\Q{{\Bbb Q}}

\def\T{{\Bbb T}}
\def\Z{{\Bbb Z}}

\def\CO{{\cal O}}
\def\CG{{\cal G}}

\def\rI{{\rm I}}
\def\rII{{\rm II}}
\def\rIII{{\rm III}}
\def\rIV{{\rm IV}}

\begin{document}

\

\vspace{1in}

\centerline{\bf  RATIONAL POINTS ON QUARTICS}

\

\centerline{\today}

\vspace{1in}

\noindent {\bf Joe Harris}

\noindent Mathematics department, Harvard University,

\noindent 1 Oxford st., Cambridge MA 02138, USA

\noindent harris{\char'100}math.harvard.edu

\

\

\noindent {\bf Yuri Tschinkel}

\noindent Mathematics department, University of Illinois at Chicago,

\noindent Chicago, IL 60607, USA

\noindent yuri{\char'100}math.uic.edu

\newpage

\

\
 
\tableofcontents

\

\

\section{Introduction}

Of all the possible extensions to higher 
dimensions of Falting's theorem, probably the most
fundamental is the 

\begin{conj}[The Weak Lang Conjecture] 
Let $X$ be a variety defined over a number field $K$.
If $X$ is of general type, then the set $X(K)$ of 
$K$-rational points of $X$ is not Zariski dense.
\end{conj}

(While the name ``Weak Lang Conjecture" has 
become standard usage---in part to distinguish
it from the ``Strong Lang Conjecture" below---we 
should point out that as stated here it was
first ventured by Bombieri (see for example \cite{noguchi}) 
and by Vojta in \cite{vojta}.)

\

We ask now whether a converse to this statement might hold. 
As it stands, the converse to the Weak Lang
Conjecture cannot possibly be true. For example, 
if we take the product $X = \P^1 \times C$ of a rational curve
and a curve
$C$ of genus $g
\ge 2$, we get a surface that is not of general type; 
but by Falting's theorem the rational points of
$X$ must lie in a finite union of fibers of $X$ over $C$.

The point is that the Kodaira dimension of a variety 
is not a sufficiently sensitive measure of
the positivity or negativity of its canonical bundle. 
If we hope to have a plausible converse to
the Weak Lang Conjecture, we need to make a stronger 
hypothesis: for example, we may make
the

\begin{conj}\label{weakconverse}
Let $X$ be a smooth projective variety defined over a number field $K$.  If the
canonical bundle $K_X$ of $X$ is negative (that is, $-K_X$ is ample), 
then for some finite
extension $K'$ of $K$ the set $X(K')$ of $K'$-rational points of 
$X$ is Zariski dense.
\end{conj}

This conjecture is easily seen to be true 
for curves and surfaces, where the hypothesis ensures
that $X$ is rational. The first real test cases are thus Fano threefolds.
In this paper, we will examine the available evidence for this conjecture, and
add to it by analyzing one further class of 
Fano threefolds, the smooth quartic hypersurfaces in
$\P^4$. Specifically, we will prove the

\begin{thm}\label{quarticthreefolds}
Let $S \subset \P^n$ be a smooth quartic 
hypersurface defined over a number field
$K$. If  $n \ge 4$, then  for some finite
extension $K'$ of $K$ the set $S(K')$ of 
$K'$-rational points of $S$ is Zariski dense.
\end{thm}

\

It is worth mentioning that Conjecture~\ref{weakconverse} 
is not the strongest possible converse to the Weak
Lang Conjecture. It may well be that we don't need 
the canonical bundle to be negative, but only nonpositive in a
suitable sense. Thus, for example, we could make the stronger

\begin{conj}\label{strongconverse}
Let $X$ be a smooth projective variety defined over a number field $K$.  If the
anticanonical bundle $-K_X$ of $X$ is nef (that is, $-K_X$ 
has nonnegative degree on every curve $C \subset X$),
then for some finite extension $K'$ of $K$ the set $X(K')$ of 
$K'$-rational points of $X$ is Zariski dense.
\end{conj}

The first test case for this conjecture would be K3 surfaces: 
it's true for curves, and looking over
the classification theorem for algebraic surfaces we may 
verify it for all surfaces except K3's. In fact, we will
prove it below for a large class of K3 surfaces, but the 
question remains open for general K3 surfaces defined
over a number field.

\

We should also mention here the Strong Lang Conjecture---that is, the

\begin{conj}[The Strong Lang Conjecture] 
Let $X$ be a variety defined over a number field $K$.
If $X$ is of general type, then there exists 
a proper subvariety $\Sigma \subsetneqq X$
such that for any finite extension $K'$ of $K$, 
$$
\#(X\setminus\Sigma)(K') \; < \; \infty \, ;
$$
that is, the set
of
$K'$-rational points of
$X$ lying outside of $\Sigma$ is finite.
\end{conj}

The converse to this statement seems plausible. 
For one thing, if $\phi : A \to X$ is any nonconstant
map from a rational or abelian variety to $X$, 
the image of $\phi$ will be defined and have a dense
collection of rational points over some extension $K'$ of $K$, 
and so will have to be contained in the
Langian exceptional subvariety $\Sigma$. Moreover, 
as a consequence of the theorem of
Kollar--Miyaoka--Mori  varieties which are not of 
general type can be thought of as  triple fibrations 
$X\ra Y\ra Z$, where   the generic fiber of $ X\ra Y$ is rationally connected, 
$Y\ra Z$ is a Kodaira fibration with generic fiber  
of Kodaira dimension $0$ and the base has  Kodaira
dimension $\le 0$ (cf. \cite{kollar-miyaoka-mori-92-2}). 
Thus the converse to the Strong Lang
Conjecture hinges on whether a variety of Kodaira 
dimension $0$ possesses a dense collection of
images of rational and/or abelian varieties. Since 
this is known for curves and surfaces, the converse
to the Strong Lang Conjecture is likewise known for 
all curvers and surfaces, and for all threefolds
except for Calabi-Yau threefolds, which represent 
the first real test. In general, however, it remains
very much open.

\

It is also worth mentioning that there are conjectures describing asymptotics for 
the number of rational points of bounded height. For example, let us consider a
smooth quartic hypersurface $S\subset \P^4$. 
Then it is expected that the number of rational points contained
in some appropriate Zariski open subset and defined over a sufficiently
large finite extension of the groundfield 
of bounded  height (induced from a standard height on $\P^4$)
grows linearly with the height (cf. \cite{FMT}).

\section{Fano Threefolds}

\

In this section we give a brief survey of 
known classification and rationality 
results for Fano threefolds (cf. 
\cite{Mori-Mukai-1}, \cite{Mori-Mukai-2}, \cite{beauville-3}, \cite{Manin}). 
It suffices to consider {\em minimal} Fano threefolds 
(not isomorphic to a blow up of a Fano variety along a smooth irreducible curve). 

\bigskip

The main invariants of Fano threefolds are: 
$r(X)$ - the index, defined as the maximal $r\in \Z$ 
such that $-K_X=rL$ for some 
$L$ in the Picard group $ \pic(X)$; $\rho(X)$ - the rank of $\pic(X)$
and the normalized degree $\delta(X)=(-K_X)^3/r(X)^3$.

\bigskip

\subsection*{\bf Group $\rI$: $r\ge 2, \rho=1$}

\medskip
\begin{enumerate}
\item $\P^3$ ($r=4$)
\item $Q^3$ ($r=3$)

The remaining 5 families have $r=2$. They are indexed by 
$\delta$. 
Let $H$ be a line bundle such that $|2H|=|-K_V|$. 

\item 
$$
\phi_H\,:\, V_1\ra \P^2
$$
is a rational map with one indeterminacy point and irreducible elliptic fibers. 
$V_1$ can be realized as a double cover of the 
Veronese cone in  $ \P^6$ whose
branch locus is a smooth intersection of this cone and a cubic hypersuface not
passing through the vertex of the cone. Another realization 
is as a hypersurface of degree 6 in the weighted projective space
$\P(1,1,1,2,3)$. 
General $V_1$ are non-rational.  Unirationality is unknown. 
\item 
$$
\phi_H\,:\, V_2\ra \P^3
$$
is a double covering ramified along a smooth quartic hypersurface.
All are unirational, the general $V_2$ is non-rational.  
\item 
$$
\phi_H\,:\, V_3\ra \P^4
$$
is a smooth cubic hypersurface. All are unirational, all are  
non-rational. 
\item
$$
\phi_H\,:\, V_4\ra \P^5
$$
is a smooth intersection of two quadrics. All are rational. 
\item 
$V_5$ is birational to a smooth quadric $Q_3$.
\end{enumerate}

\bigskip

\subsection*{\bf Group $\rII$: $r=1, \rho=1$}

\begin{enumerate}
\item 
$$
\phi_{-K_V}\,:\, W_2\ra \P^3
$$
double solid ramified along a sextic surface. They are non-rational, unirationality
is unknown. 
\item 
$$
\phi_{-K_V}\,:\, W_4\ra \P^4
$$
smooth quartic. All are non-rational,
{\em some} are unirational. In general, unirationality 
is unknown. 
\item 
$$
\phi_{-K_V}\,:\, W_6\ra \P^5
$$
smooth complete intersection of a quadric and a cubic. All are
unirational and non-rational. 
\item
$$
\phi_{-K_V}\,:\, W_8\ra \P^6
$$
smooth complete intersection of 3 quadrics. All are unirational and 
non-rational. 
\item $W_{10} $ are all unirational, rationality is unknown: the general is
non-rational.
\item $W_{12}$, $W_{16}$, $W_{18}$, $W_{22}$ are all rational.
\item $W_{14}$ are birational to a cubic threefold: all unirational, all
non-rational.
\end{enumerate}

\bigskip

\subsection*{\bf Group $\rIII$: $\rho=2,3$}

\begin{thm} $($ \cite{Mori-Mukai-1}, p. 104 $)$
If $\rho(X)=3$, then $X$ is a conic bundle over $\P^1\times \P^1$ and
has either a divisor $D\simeq \P^1\times \P^1$ or another conic bundle
structure over $\P^1\times \P^1$. 
\end{thm}

There are 4 types of minimal Fano threefolds with $\rho=3$:
\begin{enumerate}
\item Double cover of $\P^1\times \P^1 \times \P^1$ 
ramified in a $(2,2,2)$-divisor.
\item Smooth member of 
$|L^{\otimes 2}\otimes_{{\cal O}_{\P^1\times\P^1}}{\cal O}(2,3)|$
on $\P_{\P^1\times \P^1}({\cal O}\oplus {\cal O}(-1,-1)^{\oplus 2})$
such that $X\cap Y$ is irreducible.
\item $\P^1\times \P^1\times \P^1$.
\item $\P_{\P^1\times \P^1}({\cal O}\oplus {\cal O}(1,1))$.

The remaining varieties in this group have Picard number
$\rho=2$. All are conic bundles over $\P^2$.

\item Double cover of $\P^3_*$ (=  $\P^3$ blown up in one point), 
ramified in a divisor in $|-K_{\P^3_*}|$
\item Double cover of $\P^1\times \P^2$ ramified in a $(2,2)$-divisor
\item Double cover of $\P^1\times \P^2$ ramified in a $(2,4)$-divisor
\item Hypersurface of bidegree $(2,2)$ in $\P^2\times \P^2$
\item Hypersurface of bidegree $(1,2)$ in $\P^2\times \P^2$
\item Hypersurface of bidegree $(1,1)$ in $\P^2\times \P^2$
\item $\P^1\times \P^2$
\item $\P_{\P^2}({\cal O}\oplus {\cal O}(2))$
\item $\P_{\P^2}({\cal O}\oplus {\cal O}(1))$
\end{enumerate}

\

There are many more forms of Fano varieties over
nonclosed fields. The main result
of this paper together with the 
above classification implies: If $X$ is a 
Fano variety defined over some number field $K$ 
which over $\C$ is {\em not} isomorphic to 
(a blow up of) $V_1$ or $W_2$, then 
there exists a finite extension $K'/K$ such that the set
$X(K')$ is Zariski dense. 

\

One may ask for minimal conditions which insure unirationality.  
Already for Del Pezzo surfaces of degree $1$ and $2$ it is unknown 
whether or not rational points are Zariski dense as soon as there
is at least one rational point. 

\

\section{Quartic surfaces}

In this section we will prove our 
main results: to begin with, the basic

\begin{thm}\label{basictheorem}
Let $S \subset \P^3$ be a smooth quartic 
surface defined over a number field $K$, and $L
\subset S$ a line in $\P^3$ contained in 
$S$, likewise defined over $K$. Then 

a. for some finite
extension $K'$ of $K$ the set $S(K')$ of 
$K'$-rational points of $S$ is Zariski dense; and

b. if we assume further that $L$ does not meet six or more 
other lines contained in $S$, then in fact the set
$S(K)$ of $K$-rational points of $S$ is Zariski dense.
\end{thm}

As we will see in section~\ref{higherdim}, 
Theorem~\ref{quarticthreefolds} will follow as a straightforward
corollary of this statement.

\

Our proof of Theorem~\ref{basictheorem} is based 
on an analysis of the fibration of $S$ over $\P^1$ given by
projection from the line
$L$, and of the trisection of $S \to \P^1$ given 
by the points of $L$ itself. To set up the notation,
let $S \subset \P^3$ be a smooth quartic surface and 
$L \subset S$ a line. For any plane $H \subset
\P^3$ containing $L$, let
$C_H
\subset H$ be the cubic residual to $L$ in the 
intersection of $H$ with $S$---that is, such that
$S\cdot H = L + C_H$ as divisors on $H$---and 
let $D_H = C_H \cap L$ be the intersection of $C_H$ with the
line $L$. Note that projection from the line 
$L \subset \P^3$ gives a regular map $\pi$ from the surface $S$
to the line $M \cong \P^1$ parametrizing planes through $L$, 
and that the curves $C_H$ are simply the
fibers of this map. 

Note also that for any point $p \in L$, the point $p$ will 
lie in $C_H$ if and only if $H$ is the tangent plane to $S$
at $p$. Thus the
restriction of the map $\pi : S \to M$ to $L \subset
S$ is simply the restriction to $L$ of the Gauss map on $S$, 
mapping $L$ onto the line in $(\P^3)^*$ dual to $L$.
Since the general plane $H$ containing $L$ is tangent to $S$ at the three
points of $D_H$, this map has degree $3$.  The divisors $D_H$ are
the fibers of the restriction of this map, and so in particular
\emph{the divisors
$D_H$ form a linear system on $L$ of degree $3$}. 

Similarly, for any point
$p
\in L$, let $\T_pS \subset \P^3$ be the tangent plane to 
$S$ at $p$, and let $C_p = C_{\T_pS} \subset
\T_pS$ be the cubic residual to $L$ in the intersection 
of $\T_pS$ with $S$. Let $D_p = D_{\T_pS} = C_p \cap L$ be
the intersection with $L$; note that $p \in D_p$ tautologously.

\

The key point in our argument
has to do with the relation (or lack thereof) between 
the points $p$ of intersection of the curves $C_H$ with $L$
and the hyperplane class in
$\pic(C_H)$: the basic result, which we will establish 
subject to various hypotheses in the following sections,
states that for a general $H$ and any point 
$p \in C_H \cap L$, the classes of $p$ and the line bundle
$\CO_{C_p}(1)$ are linearly independent in
$\pic(C_p)$; that is, no multiple of the point $p$ is 
linearly equivalent to any multiple of the hyperplane class on
$C_H$.

\

\section{A Chow ring calculation}\label{Chow}

We'll begin by establishing a weak form of our basic result: we'll show that if $L$
meets no other line of $S$, then for general 
$H \supset L$ the points of $C_H \cap L$ are not rationally related to
$\CO_{C_H}(1)$ in $\pic(C_H)$. The proof is a 
relatively elementary argument using a calculation in the
Neron-Severi group of an associated surface. 
In the following section, we will give a more refined analysis, which
will allow us to conclude the same statement 
subject only to the weaker hypothesis that $L$ does not meet six or
more lines of $S$; while the present argument 
will be superceded by that one, the argument here is useful for its
(relative) simplicity and applications to 
similar situations. (We will see some of these in
section~\ref{otherelliptic}.)

\

Keeping the notations of the preceding section, 
the key issue in this argument, as in the following sections, has to
do with the relation (or lack thereof) between 
the points $p$ of intersection of the curves $C_H$ with $L$ and the
hyperplane class in
$\pic(C_H)$. The basic result is the following:

\

\begin{thm}\label{chowapproach} Assume that no other lines lying on $S$ meet $L$. For a general point $p
\in L$, the classes of $p$ and the line bundle $\CO_{C_p}(1)$ are linearly independent in
$\pic(C_p)$; that is, for any positive integer $n$ we have
$$
3n \cdot p \; \not\sim \; \CO_{C_p}(n) \, .
$$
\end{thm}

\

\ni {\bf Proof}. We begin by introducing a 
basic surface associated to this configuration: the
incidence correspondence
$$
T \; = \; \bigl\{(p,q) : q \in C_p \bigr\} \; \subset \; L \times S \, .
$$
To see $T$ more clearly, note first that projection from the line $L$ gives a regular map $\phi : S
\to M$ of
$S$ to the line $M \cong \P^1$ parametrizing the pencil of planes containing $L$;
the curves $C_H$ are the fibers of this map. Similarly, the divisors cut on $L$ by the curves
$C_H$ form a base-point-free pencil of degree $3$ on
$L$; the restriction $\phi' = \phi |_L$ is the 
map associated to this pencil and correspondingly has
degree $3$. In these terms, the surface $T$ is simply the fiber product
$$
T \; = \; L \times_M S \, ;
$$
in other words, $T \to L$ is the fibration 
obtained by applying the base change $L \to M$  to the 
 fibration
$S
\to M$. In particular,
$T$ is a three-sheeted cover of $S$, branched over the union of the fibers $C_H$ of $S \to M$ such
that $C_H$ is tangent to $L$. Note that for a general pair $(S,L)$ the surface $T$ will be smooth, but
it need not be always: it will be singular exactly when some curve $C_H$ is simultaneously singular
and tangent to $L$. At worst, however, it will have isolated singularities, since by the hypothesis
that $L$ meets no other lines lying on $S$ no curve
$C_H$ can have a multiple component.

\

Note that $T \to L$ has a tautologous section
$$
\Sigma \; = \; \bigl\{(p,p) : p \in L \bigr\} \; \subset \; T \, ;
$$
this is just the intersection of $T = L \times_M S \subset S \times_M S$ with the diagonal $\Delta
\subset S \times_M S$. As a divisor, the pullback $\nu^{*}(L)$ of the line $L$ under the
three-sheeted covering $\nu : T \to S$ is thus a sum
$$
\nu^{*}(L) \; = \; \Sigma + R
$$
with $R \subset T$ flat of degree $2$ over $L$. Note that since the divisors $D_H$ form a
base-point-free linear series, the general divisor $D_H$ is reduced; in particular, $R$ does not
contain $\Sigma$. $R$ will in general be irreducible, but need not be always: it will be reducible
exactly when the covering $L \to M$ is cyclic.

\

Now, suppose that the conclusion of our Theorem is false: that is, for some $n$, we have 
$$
3n \cdot p \; \sim \; \CO_{C_p}(n) 
$$
for general $p \in L$. Fixing a plane $\Gamma \subset \P^3$, there is thus for general $p$ a
rational function on
$C_p$ with a pole of order $3n$ at $p$ and zeroes of order 
$n$ at the points of intersection
$\Gamma \cap C_p$, and which is nonzero and regular everywhere else. 
It follows in turn that
there is a rational function $f$ on $T$ with divisor
$$
(f) \; = \; -3n\cdot \Sigma + n \cdot \nu^*\Gamma + D
$$
where $D$ is supported on a finite union of fibers of $T \to L$. 
Since the hypothesis that $L$ meets
no other line of $S$ insures that all fibers of $T \to L$ are 
irreducible, $D$ must consist of a sum of
fibers $C_p$ of $T \to L$. Since all fibers of $T \to L$ are 
linearly equivalent, our theorem will thus
follow from the

\

\begin{lm}\label{independent} The classes $\sigma$, 
$\gamma$ and $\phi \in A_1(T)$  of the
divisors
$\Sigma$,
$\nu^*\Gamma$ and
$C$ are independent in the group $A_1(T)$ of Weil 
divisors mod linear equivalence on $T$.
\end{lm}

\

\ni {\bf Proof}.
We need to begin with a basic fact, whose proof is 
mapped out in  Fulton, Examples 7.1.16 and
8.3.11 \cite{Fulton}.

\

\begin{lm}  Let $T$ be a normal surface. We may define, for every point $p \in T$ and
Weil divisors $D, E \in Z_1(T)$ whose supports have no common component in a neighborhood of
$p$, an intersection multiplicity
$$
j(p, D\cdot E) \; \in \; \Q \, ,
$$
bilinear in $D$ and $E$, and a bilinear intersection pairing
$$
(\,\cdot\, ,\cdot\, ) \; : \; A_1(T) \times A_1(T) \; \longrightarrow \; \Q
$$
on the group $A_1(T)$ of Weil divisors on $T$ mod rational equivalence, with the following
properties:
\begin{enumerate}
\item If $D$ and $E$ are effective and both contain $p$, then $j(p, D\cdot E) > 0$.
\item If $D$ is locally principle at $p$, that is, $D = (f)$ for some rational function $f$ in a
neighborhood of $p$, and $E$ is effective and irreducible, then
$$
j(p, D\cdot E) \; = \; \ord_p(f|_E) \, .
$$
\item If $D$ and $E$ have no common components, then
$$
([D] \cdot [E]) \; = \; \sum_{p \in D \cap E} j(p, D\cdot E) \, .
$$
\end{enumerate}
\end{lm}

Note that by the second and third properties, the intersection 
pairing extends the cup product on
the group
$\pic(T)$. 

\

We may now establish Lemma~\ref{independent} by 
calculating the matrix of intersection
products of the classes $\sigma$, $\gamma$ and 
$\phi \in A_1(T)$ and showing that this matrix is
nonsingular. All but one of these numbers are readily calculated. 
To begin with, $\phi$ is the class
of a fiber of the map $T \to L$, so of course $\phi^2=0$; 
and inasmuch as $\sigma$ is the class of a
section $\Sigma$ of that map, we have $(\phi \cdot \sigma)=1$. Next,
$\gamma$ is the pullback of the hyperplane class under the map $\nu : T \to S \hookrightarrow
\P^3$; since the map $T \to S$ has degree 3, we have
$$
\gamma^2 \; = \; 3 \cdot \deg(S) \; = \; 12 \, .
$$
Since the curves $C_p$ map forward to plane cubics under the map $\nu$, moreover, we have
$$
(\gamma \cdot \phi) \; = \; 3
$$
and similarly, since the curve $\Sigma$ maps one-to-one onto the line $L \subset S$,
$$
(\gamma \cdot \sigma) \; = \; 1 \, .
$$
In sum, then, we have the following table of intersection products
\begin{center}
\renewcommand{\arraystretch}{1.25}
\begin{tabular}{ l | r | r | r | }
& $\gamma$ & $\phi$ & $\sigma$ \\
\hline
 $\gamma$ & 12 & 3 & 1 \\
\hline
$\phi$ & 3 & 0 & 1 \\
\hline
$\sigma$ & 1 & 1 & $\sigma^2$ \\
\hline
\end{tabular}
\end{center}

\

The only mystery is the self-intersection 
$\sigma^2$ of the curve $\Sigma$ on $T$. To find this,
we use a relation of linear equivalence between 
$\Sigma$ and a curve not containing $\Sigma$: as
we saw above,
$$
\nu^*L \; = \; \Sigma + R
$$
so that if $\rho = [R] \in A_1(T)$ is the class of $R$, we have
$$
\sigma^2 \; = \; (\sigma \cdot [\nu^*L] - \rho) \, .
$$
Now, by the push-pull formula,
$$
(\Sigma \cdot \nu^*L)_T \; = \; (\nu_*\Sigma \cdot L)_S  \; = \; (L \cdot L)_S  \; = \;  -2 
$$
and so
$$
\sigma^2 \; = \; -2 -(\sigma \cdot \rho) \, .
$$
Alternatively, if we choose  the plane $\Gamma \subset \P^3$ to contain $L$, we see that the
inverse image in $T$ will consist of $\Sigma$, $R$ and the three fibers of the map $T \to L$ over
the points of intersection of $C_\Gamma$ with $L$. Thus
$$
\gamma \; = \; \sigma + \rho + 3\phi 
$$
and
$$
\sigma^2 \; = \; (\sigma \cdot \gamma - 3\phi - \rho) \; = \; -2 -(\sigma \cdot \rho) \, . 
$$

At this point, we may readily complete the calculation for general $S$ and $L$: in general,
$T$ will be smooth, and $\Sigma$ and $R$ will intersect transversely over the points $p$ where
$C_p$ is tangent to $L$ at $p$. The curves $C_H$ cut out a general pencil of degree
$3$ on
$L$, which by Riemann-Hurwitz will have $4$ branch points; thus $(\Sigma \cdot R) = 4$.
For arbitrary $S$ and $L$, however, $T$ may be singular at the points of intersection of $\Sigma$
with $R$---it will be so exactly when a curve $C_H$ has a singularity at a point of $L$---and we
can no longer say precisely what the intersection multiplicity is. All we do know, in fact, is that
$\Sigma$ must meet $R$ somewhere, so that
$$
(\Sigma \cdot R) \; > \; 0
$$
and correspondingly
$$
\sigma^2 \; < \; -2 \, .
$$

Now, we may calculate the determinant of the matrix of pairwise intersection of the classes
$\sigma$, $\gamma$ and $\phi$: it is
$$
-12 + 3 -9\sigma^2 + 3 \; = \; -6 -9\sigma^2 \; > \; 0 \, .
$$
The matrix is thus nonsingular, 
the classes $\sigma$, $\gamma$ and $\phi$ are independent, and
Lemma~\ref{independent} is proved.

\

\section{The argument via monodromy}

A finer analysis of the fibration 
$\pi : S \to M \cong \P^1$, and specifically of the monodromy of the family on
the torsion points in the Jacobian of the general fiber, yields a stronger result.

\subsection{Generalities about elliptic K3's}

To begin, we recall some basic facts about the 
fibration $S \to M$. First, since $S$ is smooth, the Gauss map $\CG :
S
\to (\P^3)^*$ is regular, and hence is finite. 
It follows that no plane $H \subset \P^3$ can be tangent to $S$ along a
curve; in other words, every hyperplane section of $S$ is reduced. 
The same is thus true of the fibers $C_H$ of
the fibration $\pi : S \to M$.

We may thus list the possible singular fibers of $\pi$: they are

\

\ni $\bullet$ \thinspace a cubic with one, two or three nodes: 
that is, an irreducible nodal curve, the union of a
line and a conic meeting transversely, or the union of three 
nonconcurrent lines. These are called \textit{fibers of
type}
$\rI_b$ with $b = 1$, $2$ and $3$ respectively;

\

\ni $\bullet$ \thinspace a cuspidal cubic, called a \textit{fiber of type} $\rII$;

\

\ni $\bullet$ \thinspace the union of a line and a tangent conic, 
called a \textit{fiber of type} $\rIII$; or

\

\ni $\bullet$ \thinspace the union of three concurrent lines, 
called a \textit{fiber of type} $\rIV$.

\

Note that fibers of type $\rI_b$ correspond to poles of order $b$ 
of the $j$-function on $M$ associated to the
elliptic fibration $S \to M$. By contrast, fibers of type $\rII$ 
and $\rIV$ correspond to zeroes of $j$---after base
changes of order $6$ and $3$ respectively, we may replace the 
singular curve by an elliptic curve  $\tilde C$ of
$j$-invariant $0$---and fibers of type $\rIII$ correspond to 
points where $j = 1728$: after a base change of
order $4$ we may replace the singular curve by an elliptic curve $\tilde C$ of
$j$-invariant $1728$. The monodromy action on the homology 
of the general fiber around a fiber of type $\rI_b$
is thus given by the Picard-Lefschetz transformation, while 
in the case of fibers of type $\rII$, $\rIII$ and
$\rIV$ the monodromy is the just the action of 
the automorphism given by the base change on $\tilde C$. We list
here these actions; see the original paper of Kodaira  
or the discussion in Barth--Peters--Van de Ven  
for details (cf. \cite{barth-peters-vdv}). 

\

\begin{center}
\begin{tabular}{ c | c }
type & monodromy \\
\hline
 $\rI_b$ & \parbox[c][.6in]{1.3in}{$M_{\rI_b} = \begin{pmatrix}
1 & b \\
0 & 1
\end{pmatrix}$} \\
\hline
$\rII$ & \parbox[c][.6in]{1.3in}{$M_{\rII} = \begin{pmatrix}
1 & 1 \\
-1 & 0
\end{pmatrix}$} \\
\hline
$\rIII$ & \parbox[c][.6in]{1.3in}{$M_{\rIII} = \begin{pmatrix}
0 & 1 \\
-1 & 0
\end{pmatrix}$} \\
\hline
$\rIV$ & \parbox[c][.6in]{1.3in}{$M_{\rIV} = \begin{pmatrix}
0 & 1 \\
-1 & -1
\end{pmatrix}$} \\
\end{tabular}
\end{center}

There is one constraint on the number of singular fibers $C_H$: 
by a standard Euler characteristic calculation,
$$
\chi(S) \; = \; 24 \; = \; \sum_{H \in M} \chi(C_H)\, .
$$
The Euler characteristics of fibers of type $\rI_b$, $\rII$, 
$\rIII$ and $\rIV$ are $b$, $2$, $3$ and $4$
respectively, giving a linear relation on the numbers of 
fibers of each type. In particular, we see that there must
be at least $6$ singular fibers $C_H$.

\

A related issue is the count of curves $C_H$ that are not 
transverse to $L$. Given that the divisors $D_H$ cut
out on $L$ by the curves $C_H$ form a pencil of degree $3$, 
it follows by Riemann-Hurwitz that there must be a
total of $4$ branch points, counting multiplicity: 
that is, either two curves $C_H$ having a point of intersection
multiplicity $3$ with $L$, one such curve and two others 
having a double point of intersection with $L$, or four
curves $C_H$ having a double point of intersection with $L$. 
Moreover, note that the monodromy on the points
$p_i$ around a curve $C_{H_0}$ having a point $p$ 
of intersection multiplicity $2$ (respectively, $3$) with $L$ is
necessarily cyclic of order exactly $2$  (respectively, $3$). 
It follows that if the points $p$ of $C_H \cap L$ differ
from each other by torsion in $\pic(C_H)$ for general $H$, 
then since sections of an elliptic fibration that differ by
torsion in the general fiber can intersect only at singular points of fibers,
\textit{the multiple point
$p$ of intersection of
$C_{H_0}$ with $L$ must be a singular point of $C_{H_0}$}.

On the other hand, since there can be at most four fibers 
$C_H$ having multiple points of intersection with
$L$, and each can contribute at most $4$ to the Euler 
characteristic,  we may draw one conclusion in particular
that will turn out to be vital to the following analysis: that
\textit{there must be singular fibers
$C_H$ that intersect
$L$ transversely}.

\

\subsection{Analysis of the points of $C_H \cap L$}

We now ask the crucial question: is it possible 
that for general planes $H \supset L$ containing $L$, the points of
intersection of $C_H$ with $L$ are all rationally 
related to the class $\CO_{C_H}(1) \in \pic(C_H)$? For the
remainder of this subsection, then, we will 
assume that this is the case, and see what conclusions we may draw
from it.

So: let $H$ be a general plane containing $L$, and let 
$p_1$, $p_2$ and $p_3$ be the three points on intersection
of
$C_H$ with $L$. Assume that some multiple of $p_i$ is 
linearly equivalent to a multiple of the hyperplane section
of $C_H$, and let $n$ be the smallest positive integer such that
$$
3n \cdot p_i \; \sim \CO_{C_H}(n) \, .
$$
Note that since the monodromy on the three points $p_i$ 
as $H$ varies is at least transitive, this hypothesis will
hold for one $p_i$ if and only if it holds for all three, 
and the value of $n$ will be the same for all three. In this
case, the pairwise differences
$$
\alpha_{i,j} \; = \; p_i - p_j \; \in \; \pic^0(C_H)
$$
are torsion; we let $m$ be their order (again, 
since the monodromy is transitive on the three pairs
$\{\pm\alpha_{i,j}\} = \{\alpha_{i,j}, \alpha_{j,i}\}$, 
the value of $m$ will be the same for all $i$ and $j$). Note
that the classes $\alpha_{i,j}$ are all nonzero, but 
need not be distinct: if $m=2$, of course, we have $\alpha_{i,j} =
\alpha_{j,i}$, while if $m=3$ we could have 
$\alpha_{1,2} = \alpha_{2,3} = \alpha_{3,1}$. If $m > 3$, however, we
can see from the transitivity of the monodromy and the 
fact that $\alpha_{1,2} + \alpha_{2,3} + \alpha_{3,1} = 0$
that they must be distinct. 

Note finally that if $m \neq 2$, the monodromy on 
the points $p_i$ must be cyclic,
rather than the symmetric group $S_3$: a transformation 
fixing $p_1$, for example, and exchanging $p_2$ and
$p_3$, would exchange $\alpha_{1,2}$ and  $\alpha_{1,3}$ 
and send  $\alpha_{2,3}$ to  $-\alpha_{2,3}$; given
that $\alpha_{1,2} - \alpha_{1,3} +\alpha_{2,3}=0$, 
this implies $2\alpha_{2,3}=0$. It follows in particular that in
case $m \neq 2$, there are exactly two fibers $C_H$ 
not transverse to $L$, and each has a triple point of
intersection with $L$.

\

To carry out the further analysis of the monodromy 
action on the points $p_i$ and the classes $\alpha_{i,j}$, we
will consider in turn three potential cases: $m > 3$, $m=3$ and $m=2$.

\

\ni \underbar{Case a: $m > 3$}.  The crucial 
observation here is simple enough: it is that if $C_H$ is any
fiber transverse to $L$, then the monodromy 
around $C_H$ acts trivially on the points $p_i$ and hence on the
classes
$\alpha_{i,j}$. Let us consider accordingly the 
fixed point sets of the various monodromy transformations listed
above:

\begin{enumerate}
\item The transformation $M_{\rI_b}$ visibly fixes a 
cyclic subgroup of the points of order $m$ in $\pic^0(C_H)$,
and fixes no other primitive point of order $m$.
\item $M_\rII$ has characteristic polynomial 
$p(\lambda) = \l^2 - \l +1$, which does not
vanish at $1$; thus $M_\rII$ does not fix any primitive points of order $m$.
\item $M_\rIII$ has characteristic polynomial 
$p(\lambda) = \l^2  +1$, and since we have assumed $m \neq 2$
this does not vanish at $1$; thus $M_\rIII$ does not fix any primitive points of order $m$.
\item $M_\rIV$ has characteristic polynomial 
$p(\lambda) = \l^2 + \l +1$, and since we have assumed $m \neq 3$
this does not
vanish at $1$; thus $M_\rIV$ does not fix any primitive points of order $m$.
\end{enumerate}

The conclusion is plain: all singular fibers $C_H$ 
transverse to $L$ must be of type $\rI_b$; and since we have
already seen that there must exist such fibers, 
we deduce further that \textit{the classes $\alpha_{i,j}$ lie in a
cyclic subgroup $G \cong \Z/m$ of $\pic^0(C_H)$}.

\

We ask now what may be the singular fibers $C_H$ 
that are singular at a point of $L$. As we have seen,
since $m \neq 2$ there are exactly two fibers 
$C_H$ not transverse to $L$, and each has a triple point of
intersection with
$L$. For each such fiber, moreover, the subgroup 
$G \subset \pic^0(C_H)$ spanned by (any one of) the
$\alpha_{i,j}$ is an eigenspace for the action of the 
monodromy on the points of order $m$ in $\pic^0(C_H)$. Now,
the monodromy action $M_{\rI_b}$ associated to a singular 
fiber of type $\rI_b$ has only one eigenspace, with
eigenvalue $1$. But the monodromy on the points $p_i$ is 
necessarily nontrivial and cyclic, and since $m\neq 3$
it follows that the action on the classes $\alpha_{i,j}$ is nontrivial. 
Thus a fiber of type $\rI_b$ singular at a point
of $L$ cannot occur.

Next,  suppose  that we have a
singular fiber of type $\rII$. As we saw, the monodromy $M_\rII$ 
has characteristic polynomial $p(\l) = \l^2 - \l
+ 1$. Suppose $\l$ is a root of this polynomial mod $(m)$. Then we have
$$
\l^2 \; \equiv \; \l - 1 \quad {\rm mod} \; m
$$
and hence
$$
\l^3 \; \equiv \; \l^2 - \l  \; \equiv \; -1 \quad {\rm mod} \; m \, .
$$
Alternatively, we can just multiply out and 
see that $M_\rII^3 = -1$. Either way, we see that $M_\rII^3$ cannot
fix any element of $\pic^0(C_H)$ of order $m \neq 2$, and so no fiber of this type can occur.

A similar analysis shows that no fiber of 
type $\rIII$ can occur (or we could just observe that the union of a line
and a tangent conic cannot have a point of 
intersection multiplicity $3$ with $L$); and so we conclude in sum that
\textit{every singular fiber $C_H$ transverse to 
$L$ must be of type $\rI_b$; and every fiber not transverse to $L$
must be of type $\rIV$}. In particular, 
we see that in this case $L$ must meet at least the $6$ other lines of $S$
comprising the two fibers of type $\rIV$.

\

We do not know of any examples of this case.

\

\ni \underbar{Case b: $m = 3$}.  The analysis here 
is very similar to the preceding case, with one
exception. To begin with, the same analysis shows 
that the fibers $C_H$ not transverse to $L$
must both be of type $\rIV$: none of other 
transformations $M_{\rI_b}$, $M_\rII$ and $M_\rIII$ has a
three-cycle on the points of order $3$. One difference 
is that it is a priori possible that the classes $\alpha_{i,j}$
do not lie in a cyclic subgroup: this could be the case 
if every singular fiber $C_H$ transverse to $L$ is of type
$\rI_3$. But if there were $\delta$ such fibers, the 
formula for the Euler characteristic gives
$$
24 \; = \; 2\cdot 4 + \delta \cdot 3
$$
which has no solution. Thus the classes $\alpha_{i,j}$
 lie in a cyclic subgroup, which means in turn 
that we must have $\alpha_{1,2} = \alpha_{2,3} = \alpha_{3,1}$.
Thus in turn implies that
$$
3p_i \; \sim \; p_1 + p_2 + p_3 \; \sim \; \CO_{C_H}(1) \, ;
$$
that is, $n=1$, or in other words 
\textit{all three points $p_i$ are flexes of $C_H$} for general $H$. Again, we do not know
if this is possible. In any event, we see again 
that this case can occur only when six or more other lines of $S$
meet $L$.

\

\ni \underbar{Case c: $m = 2$}.  This case 
gives us the least amount of control over the behavior of the
singular fibers
$C_H$ not transverse to $L$, but happily the 
most over the behavior of those that are. Very simply, in this case 
the classes $\alpha_{i,j}$ cannot
 lie in a cyclic subgroup: they must comprise 
all three classes of order $2$ in $\pic^0(C_H)$. It follows in turn that
the monodromy associated to each fiber $C_H$ 
transverse to $L$ must be trivial on  the points of order $2$ in
$\pic^0(C_H)$, which says in turn that 
\textit{every singular fiber  $C_H$ transverse to $L$ must
be of type $\rI_2$}, that is, must consist of 
the union of a line and a conic. Now, as far as we can tell, the fibers
$C_H$ not transverse to $L$ can be of type 
$\rII$, $\rIII$ or $\rIV$, but in any event the total contribution of
such fibers to the Euler characteristic 
can be at most $8$ (if we have either four fibers of type $\rIII$ or
two of type $\rIV$). The remaining singular fibers 
must therefore contribute at least $16$ to the Euler
characteristic, which means we must have at least 
$8$ fibers of type $\rI_2$; in particular, \textit{$L$ must
intersect at least $8$ other lines of $S$}.

\

In sum, we have established the

\begin{thm}\label{geometriclemma}
Let $S \subset \P^3$ be a smooth quartic surface  and $L
\subset S$ a line in $\P^3$ contained in $S$; 
assume that $L$ does not meet six or more
other lines contained in $S$. Let $n$ 
be any positive integer. If $p \in L$ is a general point, then
$$
3n\cdot p \; \not\sim \; \CO_{C_p}(n) \, .
$$
\end{thm}

\

\section{Rational points on quartic surfaces}

To deduce part b) of Theorem~\ref{basictheorem} 
from the analysis carried out so far, we need one more
ingredient. Briefly, Theorems~\ref{chowapproach} and \ref{geometriclemma}
assure us (subject to their hypotheses)
that for a general point $p \in L(\overline \Q)$, the cubic $C_p$
has positive rank, and hence a dense set of rational points, over the field
of definition of $p$. It seems reasonable to expect that for ``most" of the 
points $p \in L(K)$ this would be true. But there will be countably many
points $p \in L(\overline \Q)$ for which $p$ \underbar{is} rationally related
to $\CO_{C_p}(1)$---for each $n$, there is a finite 
subset $\Phi_n \subset L$ $p \in L(\overline \Q)$ such that $3n\cdot p
\sim \CO_{C_p}(n)$ in $\pic^0(C_p)$---and 
it is still a logical possibility,
if not a plausible one, that all the points of $L(K)$ lie in this set.

\

There are two ways of eliminating this possibility. The first is to invoke
an extremely powerful theorem due to Merel \cite{Merel}:

\begin{thm}[Merel]\label{merel}
Let $K$ be any number field. There is an integer 
$n_0=n_0(K)$ such that no elliptic curve defined over $K$ has a
$K$-rational point of order $n > n_0$.
\end{thm}

This theorem assures us that for $n > n_0$, the subset $\Phi_n$ is disjoint
from $L(K)$; so that for all but finitely many $p \in L(K)$ the point $p$ will not be
rationally related to $\CO_{C_p}(1)$ in $\pic(C_p)$.

\

A second way to arrive at this fact uses a theorem 
that is easier to prove, if less simple to state.
Suppose that $E \to B$ is any family of elliptic curves 
and  $\sigma$ a section of
$E \to B$, both defined over a number field $K$. For each point
$t
\in B$, let
$h_B(t)$ be the height of $t$ (relative to any 
divisor of degree $1$ on $B$), $h(\sigma_t)$ the
canonical height of the value $\sigma_t$ of the 
section in the fiber $E_t$ of the family over $t$, and
$h(\sigma)$ the canonical height of the section 
$\sigma$ in the fiber of $E$ over the generic point of
$B$. We have then the following theorem of 
Dem'janenko \cite{dem}, Manin \cite{Manin2} and
Silverman \cite{sil}:
\begin{thm}
$$
\lim_{h_B(t) \to \infty} \;  \frac{h(\sigma_t)}{h_B(t)} \;  = \;  h(\sigma) \, .
$$ 
\end{thm}

In particular, to say that $\sigma$ is not a torsion point in the fiber of $E$ over the generic point of
$B$ is to say that $h(\sigma) > 0$, and we can deduce the 

\begin{cor}
If $\sigma_t$ is not a torsion point in $E_t$ for general $t$, then there are only finitely points $t \in
B(K)$ such that $\sigma_t$ is a torsion point in $E_t$.
\end{cor}

Applying this to the family $T \to L$ of elliptic curves introduced in section~\ref{Chow} (with the
origin given by the tautologous section $\Sigma$ and the section $\sigma$ given by the divisor class
$\CO_{C_p}(1)$ in each fiber, we deduce again that (subject to the hypotheses of
Theorem~\ref{geometriclemma}) for all but finitely many $p \in L(K)$ the point $p$ will not be
rationally related to $\CO_{C_p}(1)$ in $\pic(C_p)$.

\

Now, suppose $S \subset \P^3$ is a smooth 
quartic surface defined over a number field $K$ and $L
\subset S$ a line in $\P^3$ contained in $S$, 
likewise defined over $K$; assume that $L$ does not meet six or more
other lines contained in $S$  (not necessarily 
defined over $K$). For each point $p \in L(K)$ and for each integer 
$n$ there are  unique  points $q_n$ and $r_n \in
C_p$ such that
$$
q_n + (3n-1)\cdot p \; \sim \; \CO_{C_p}(n) 
$$
and
$$
-r_n  + (3n+1)\cdot p \; \sim \; \CO_{C_p}(n) \, ,
$$
also defined over $K$. Moreover, for all but finitely many $p \in L(K)$
we have 
$$
3n\cdot p \; \not\sim \; \CO_{C_p}(n) 
$$
for every $n$, so that these points are all distinct. 
We have, accordingly, an infinite collection of $K$-rational points on
$C_p$, so that $C_p$ is contained in the Zariski 
closure of $S(K)$; since this is true of infinitely many curves $C_p$,
it follows that the Zariski closure of $S(K)$ is all of $S$.

\

As for the proof of part a) of Theorem~\ref{basictheorem}, 
given part b) this requires only one further trick, and
it's a relatively simple one. This is expressed in the

\begin{lm}\label{threelines}
Let $S \subset \P^3$ be a smooth quartic surface  and $L$, $L'$ and $L''
\subset S$ three lines in $\P^3$ contained in $S$; 
assume that $L$ does not meet either $L'$ or $L''$, but that $L'$
and
$L''$ do meet. For each plane $H \supset L$ containing $L$, 
let $q_H = C_H \cap L'$ and $r_H = C_H \cap L''$. Fix
an integer $n \ge 2$. If
$H \supset L$ is a general plane containing $L$, then
the difference $q_H - r_H$ is not torsion in $\pic^0(C_H)$.
\end{lm}

\ni {\bf Proof}. This is easy. Note first that since $L$ 
doesn't meet $L'$ and $L''$, no plane $H$ containing $L$ can
be the tangent plane to $S$ at any point of $L'$ or $L''$; 
in other words, $q_H$ and $r_H$ are smooth points of the
curve
$C_H$ for all $H$. 

Now let $T$ be as in the proof of Theorem~\ref{chowapproach}. 
The assignment to each point $p \in L$ of
the points
$q_p$ and
$r_p \in C_p$  give sections of the map $T \to L$, 
which we claim do not differ by a translation of finite order in
the general fiber $C_p$. But two sections of an 
elliptic fibration that differ by torsion in the general fiber can
intersect only at singular points of fibers, and so we are done.

\

\ni {\bf Proof of part a) of Theorem~\ref{basictheorem}}. 
Now suppose we have a quartic surface $S$ and a line $L
\subset S$. If fewer than
six other lines of
$S$ meet
$L$, we may apply part b) of Theorem~\ref{basictheorem} 
to conclude that the points of $S$ rational over the
field of definition of $L$ are Zariski dense. 
Suppose conversely that {\it every line of $S$ meets at least six other
lines of $S$}. We claim in that case that $S$ 
must contain a configuration of lines as in Lemma~\ref{threelines}. To
see this, start with any line $L_0 \subset S$. 
Since no more than four lines on $S$ can pass through a single point
of $s$, we have to consider only two possibilities:

\

\underbar{Case 1}: three pairwise skew lines $L_1$, 
$L_2$ and $L_3 \subset S$ meet $L_0$. In this case, let
$M_1,\dots,M_5\subset S$ be five other lines meeting 
$L_1$ (in addition to $L_0$). If any one of them fails to
meet both $L_2$ and $L_3$, we are done: if $M_i$ fails 
to meet $L_j$, we take for our configuration $L = L_j$, $L'
= L_1$ and $L'' = M_i$. On the other hand, if all six 
lines $L_0,M_1,\dots,M_5$ meet all three lines $L_1$, $L_2$
and $L_3$, they must all lie on the unique quadric 
surface $Q \subset \P^3$ containing $L_1$, $L_2$
and $L_3$---but then we have nine lines in $Q \cap S$, 
contradicting B\'ezout (among others).

\

\underbar{Case 2}: two triples of 
concurrent lines $\{L_1, L_2,L_3\}$ and $\{M_1, M_2,M_3\}$ meet $L_0$. Even
easier: just let $N$ be any line meeting $L_1$ and 
skew to $L_0$. $N$ can't meet all three lines $M_1$, $M_2$
and $M_3$ (it would be coplanar with them and hence 
meet $L_0$); if it misses $M_i$ we take $L = M_i$,
$L'=L_1$ and
$L''=N$.

\

To complete the proof of Theorem~\ref{basictheorem}, 
suppose now that $L$, $L'$ and $L'' \subset S$ are a
configuration as in Lemma~\ref{threelines}; let $K'$ 
be any field over which $S$ and all three lines are defined.
For each plane $H$ containing $L$ and each integer
$n$ there is a unique  point
$x_n
\in C_H$ such that
$$
x_n + n\cdot q_H \; \sim \; (n+1)\cdot r_H \, ,
$$
also defined over $K'$. Moreover, by Merel's theorem 
and Lemma~\ref{threelines}, for all but finitely many $H$
these points are all distinct. We have thus for 
infinitely many $H$ an infinite collection of
$K'$-rational points on
$C_H$, and once more
it follows that the Zariski closure of $S(K')$ is all of $S$.

\

\section{An example: the Fermat quartic}

In light of the analysis above, it might seem 
unlikely that there is any quartic surface $S$ and line $L$ such that
for general planes $H$ containing $L$, the points of 
$C_H \cap L$ are all rationally related to the hyperplane
class in $\pic(C_H)$. In fact, however, it does occur: 
we will describe here the unique example we know of, the
Fermat surface. See also the analysis of Piatetski-Shapiro
and Shafarevitch in \cite{shaf}.

\

To begin with, we take $S \subset \P^3$ the quartic given by the equation
$$
X^4 - Y^4 + Z^4 - W^4 \; = \; 0
$$ 
and $L$ the line
$$
X \; = \; Y \quad {\rm and} \quad Z \; = \; W \, .
$$
Any plane containing $L$, other than the plane $Z=W$, 
can be realized as the span of $L$ and a third point of the
form $[a,-a,1,-1]$ for some scalar $a$; that is, the plane given parametrically by
$$
[U,V,T] \; \longmapsto \; [U+aT, U-aT, V+T, V-T] \, .
$$
Restricting the equation of $S$ to $H$ gives the equation
$$
(U+aT)^4- (U-aT)^4+(V+T)^4-(V-T)^4 \; = \; 8aU^3T + 8a^3UT^3 + 8V^3T + 8VT^3
$$
so the equation of $C_H$ is simply
$$
F_a(U,V,T) \; = \; aU^3 + a^3UT^2 + V^3 + VT^2 \; = \; 0 \, .
$$
The points of $C_H \cap L$ are given by the 
further equation $T=0$; that is, they are the points
$$
[U,V,T] \; = \; [1, b, 0]
$$
where $b^3 = a$. Note that the monodromy on these as $a$ varies is cyclic.

\

What are the singular fibers $C_H$ of the fibration $S \to M$? 
The plane $Z=W$ corresponding to $a = \infty$ is
certainly one, consisting of three concurrent lines 
meeting a point of $L$; the plane $X=Y$ corresponding to $a =
0$ is another singular fiber of type $\rIV$. To find 
the remaining ones we have simply to write out the partial
derivatives of $F_a(U,V,T)$ and equate them all to $0$; the equations
$$
3aU^2+a^3T^2 \; = \; 3V^2+T^2 \; = \; 2VT + 2a^3UT \; =  0
$$
together imply that $a^4 = \pm1$, and that
$$
U \; = \; \frac{1}{\sqrt{-3}}T \quad {\rm and} \quad V \; = \; \frac{a}{\sqrt{-3}}T \, .
$$
There are thus $8$ singular fibers apart from $a=0$ and $a=\infty$, 
each having two singular points (that is,
consisting of a line and conic).

\

We claim now that \textit{the three points $p_i$ of 
intersection of $C_H$ with $L$ differ by torsion of order $2$,
and each satisfies $6\cdot p_i \sim \CO_{C_H}(2)$}. 
The second statement follows from the first, given that
$p_1+p_2+p_3 \sim \CO_{C_H}(1)$, and the first is 
readily checked: we simply observe that \textit{two points $p$ and
$q$ on a plane cubic curve $C$ differ by torsion of 
order $2$ if and only if the point of intersection $\T_pC \cap
\T_qC$ of the tangent lines to $C$ at $p$ and $q$ lies on $C$}. 
Now, the equation of the tangent line to $C_H$ at a
point $[\mu,\nu,\tau]$ is
$$
(3a\mu^2+a^3\tau^2)\cdot U\; +\; 
(3\nu^2+\tau^2)\cdot V \; + \; (2\nu\tau + 2a^3\mu\tau)\cdot T \; = \; 0
$$
and at the point $[1,b,0]$ for some cube root $b$ of $a$ this is
$$
-3b^3\cdot U + 3b^2\cdot V \; = \; 0 \, .
$$
For any two distinct cube roots $b$ of $a$, 
the resulting linear forms in $U$ and $V$ are independent, so that the
point of intersection of any two of the tangent 
lines to $C_H$ at the points of $C_H \cap L$ is just
$$
[U,V,T] \; = \; [0, 0, 1]
$$
which is a point of $C_H$.

\

Again, we don't know of any other examples of 
a quartic surface $S$ and a line $L \subset S$ such that
for general planes $H$ containing $L$ the points 
of $C_H \cap L$ are all rationally related to the hyperplane
class in $\pic(C_H)$, nor do we know any examples 
at all where the points differ from each other by torsion of
order greater than $2$.

\

\section{Quartic threefolds}\label{higherdim}

We would now like to use Theorem~\ref{basictheorem} 
to deduce Theorem~\ref{quarticthreefolds}. This is
relatively simple: we just have to check that if 
$X \subset \P^n$ is any smooth quartic hypersurface and $L
\subset X$ any line, then for a general $3$-plane 
$\P^3 \subset \P^n$ containing $L$, the surface $S = X \cap
\P^3$ and the line $L$ satisfy the hypotheses of 
part b) of Theorem~\ref{basictheorem}. It's enough to do this in
case $n=4$, and it requires only a straightforward geometric argument.

\

We start with a basic fact:

\begin{lm}\label{dimone} If $X \subset \P^4$ is a 
smooth quartic hypersurface, the Fano variety
$F_1(X)
\subset
\G(1,4)$ of lines on
$X$  has pure dimension one.
\end{lm}

\ni {\bf Proof}. To begin with, the homogeneous  
quartic polynomial $F \in \Sym^5(\C^5)$
on $\P^4$ defining $X$ gives rise to a section 
$\tau_F$ of the fifth symmetric power $\Sym^5(S^*)$
of the dual of the universal subbundle $S$ on 
$\G(1,4)$, and the zero locus of this section is the
Fano scheme $F_1(X)$. This shows that it has 
dimension at least 1 everywhere, and  since the top
Chern class
$c_5(\Sym^5S^*) = 320\sigma_{3,2}
\ne 0$ that it is nonempty; it remains only 
to see that it cannot be two- or higher-dimensional.

To do this, let 
$$
\Phi = \P(S|_{F_1(X)}) = \{(L,p) : p \in L\} \subset F_1(X) \times X
$$ 
be
the univeral projective line bundle over the 
Fano variety $F_1(X) \subset \G(1,4)$, and $\rho : \Phi
\to X
\subset \P^4$ the projection map. The tangent 
space to $F_1(X)$ at a point $L \in F_1(X)$ may be
identified with the space of section $H^0(L,N_{L/X})$; 
and in these terms, at a general point $(L,p)
\in
\Phi$ the image of the differential 
$d\rho_{(L,p)} : \T_LF_1(X) \to \T_pX$ mod $\T_pL$ is simply the
image of the  map $H^0(L,N_{L/X}) \to (N_{L/X})_p = \T_pX/\T_pL$ 
given by evaluation at $p$. Since
$N_{L/X}$ has exactly one summand of 
nonnegative degree, this image is always one-dimensional
(mod $\T_pL$), and so we may conclude that 
the image of the map $\rho : \Phi \to X$---that is, the
union of the lines on $X$---is always exactly $2$-dimensional. But
$X$ contains no $2$-planes, and no surface in 
projective space other than a $2$-plane may contain
$\infty^2$ lines, so we may conclude that
$X$ contains only $\infty^1$ lines.

\

On the basis of a naive dimension count, we 
would expect the map $\rho$ from the
two-dimensional variety $\Phi$ to the threefold $X$ to 
have a one-dimensional double point locus,
that is, $\infty^1$ pairs of distinct lines $L, L' \subset X$ 
that meet; we'd thus expect that each
line of $X$ would meet finitely many others. (We'll calculate the number in just a
moment.) Accordingly, we'll call a line
$L
\subset X$ {\it exceptional} if it meets infinitely many 
other lines of $X$. We will denote by
$F_1^e(X) \subset F_1(X)$ the locus of exceptional lines.

By a straightforward dimension
count, a general quartic threefold
$X$ has no exceptional lines, and a general $X$ containing 
an exceptional line will contain only
finitely many; in particular, it will contain nonexceptional 
lines as well. Our situation is that we're
able to apply Theorem~\ref{basictheorem} directly to the 
surface $S = H \cap X$ where $H$ is a general
hyperplane section containing a nonexceptional line $L$, 
and so our concern is whether an arbitrary $X$ may
contain only exceptional lines. In fact, this is possible: 
the Fermat quartic is one example (we don't know any
other examples). What we want to do, accordingly, 
is to say as much as we can about quartic threefolds that have
positive-dimensional families of exceptional lines.

So: let $X$ be a smooth quartic threefold, and 
$\Psi \subset F_1^e(X)$ an irreducible component of
the Fano variety of lines on $X$ consisting entirely 
of exceptional lines. For each line $L \in \Psi$,
the locus of lines meeting $L$ will contain one or 
more irreducible components of $F_1(X)$, and so
we must have one of the following two situations:
\begin{enumerate}
\item There are two irreducible components 
$\Psi, \Psi' \subset F_1(X)$ such that every pair of
lines $L \in \Psi$ and $L' \in \Psi'$ meet; or
\item There is an irreducible component 
$\Psi \subset F_1(X)$ such that every pair of lines $L, L'
\in \Psi$ meet.
\end{enumerate}
In the first case, the surface $S \subset X$ swept 
out by the  lines of $\Psi$ has two rulings by
lines; but the only surface in projective space with 
two rulings by lines is a quadric surface in
$\P^3$, and since $\pic(X) = \Z\langle \CO_X(1) \rangle$, 
$X$ contains no such surfaces; thus the
first case cannot occur. In the second case, 
let $L, L' \in \Psi$ be two general lines, meeting at a
point $p$. A third general line $L'' \in \Psi$ must 
meet both $L$ and $L'$; if doesn't pass through
$p$ it must lie in the plane spanned by $L$ and $L'$, 
and so this plane would have to be contained
in $X$. Since $X$ contains to $2$-planes, we conclude 
that {\em all the lines $L \in \Psi$ have a
common point $p$}. It follows in particular that 
all the lines $L \in \Psi$ are contained in the
tangent plane $\T_pX$, and hence that $X \cap \T_pX$ 
is simply the cone with vertex $p$ over an
irreducible plane quartic curve. Finally, since $X$ 
is smooth the Gauss map $\CG$ cannot be
constant on any curve of $X$; thus $X \cap \T_pX$ 
can have at most isolated singularities and hence
the curve $C$ must be smooth. We have thus established the

\begin{prop}\label{exceptionallines}
Let $X$ be a smooth quartic threefold, and 
$\Psi  \subset F_1^e(X)$ an irreducible component of
the Fano variety of lines on $X$ consisting 
entirely of exceptional lines. Then there is a point $p \in
X$ such that $X \cap \T_pX$ is the cone with vertex 
$p$ over a smooth plane quartic curve, and
$\Psi$ is simply the ruling of this cone. In particular, $\Psi
\subset \G(1,4) \subset \P^9$ has as underlying 
reduced scheme a smooth plane quartic curve.
\end{prop}

It's a nice exercise to check directly that in 
this case the component $\Psi$ of the Fano scheme
$F_1(X)$ has multiplicity $2$. 
Since by the proof of Lemma~\ref{dimone}
the Fano scheme is a curve of degree $320$ in $\G(1,4) \subset \P^9$, 
it follows that {\em a quartic threefold $X$
will contain only exceptional lines if and only if it has 
exactly $40$ hyperplane sections consisting of
cones over quartic plane curves}. Again, this is the case 
for the Fermat quartic; we don't know if
there are others.

\

In any event, Theorem~\ref{quarticthreefolds} now 
follows readily from Theorem~\ref{basictheorem}. Let $X
\subset \P^4$ be any quartic threefold defined over a field $K$. 
Suppose first that $X$ contains a non-exceptional
line
$L$; say $L$ is defined over a field
$K' \supset K$. Then for a general hyperplane 
$H \cong \P^3 \subset \P^4$ containing $L$, the surface $S_H = X
\cap H$ will contain no other lines meeting $L$, and
by part b) of
Theorem~\ref{basictheorem} we may deduce that $S_H(K')$ 
is Zariski dense in $S_H$; hence $X(K')$ is Zariski
dense in $X$.

If on the other hand every line of $X$ is exceptional, 
then by Proposition~\ref{exceptionallines} the Fano variety
of lines on $X$ will consist of a union of curves 
supported on plane quartic curves. If $L \subset X$ is a line
corresponding to a general point $[L] \in \Psi \subset F_1(X)$ 
of a component $\Psi$ of $F_1(X)$---in particular,
if it is not a point of intersection of $\Psi$ 
with another component of $F_1(X)$---then $L$ can meet only finitely
many lines $L' \subset X$ corresponding to 
points $[L'] \notin \Psi$. It follows that for a general hyperplane $H
\cong \P^3 \subset \P^4$ containing $L$, the surface $S_H = X
\cap H$ will contain exactly three other lines 
meeting $L$; so once again the hypothesis of part b) of
Theorem~\ref{basictheorem} is satisfied, and we 
conclude that the points of $X$ rational over the field of
definition of $L$ are Zariski dense.

\

\section{Other elliptic surfaces}\label{otherelliptic}

As suggested in section~\ref{chowapproach}, the approach 
via the calculation of intersection numbers in the
Neron-Severi group of an associated surface works in 
substantially greater generality. Specifically, we will prove

\begin{thm}\label{ellipticlemma}
Let $S$ be any smooth irrational surface and $\pi : S \to \P^1$ 
an elliptic fibration all of whose fibers are 
irreducible; for $\lambda \in \P^1$, let $E_\l = \pi^{-1}(\l)$ 
be the fiber of $S \to \P^1$ over $\l$. Let $C \subset
S$ be any smooth rational or elliptic curve, of degree
$m
\ge 2$ over
$\P^1$ and $n \neq 0$ any integer. Then for general $\l \in \P^1$ and any $p \in E_\l \cap C$, 
$$
nm\cdot p \; \not\sim \; \CO_{E_\l} \otimes \CO_S(n\cdot C) \, .
$$
\end{thm}

As before, we may deduce from this immediately the

\begin{cor}\label{ellipticcor}
Let $S$ be any smooth surface, $\pi: S \to \P^1$  an 
elliptic fibration with irreducible fibers as in
Theorem~\ref{ellipticlemma}.

\ni a) Let $C \subset S$ be a smooth rational curve. Assume that
$S$,
$\pi$ and
$C$ are defined over a field $K$ and that $C$ is rational over 
$K$. Then the set $S(K)$ of $K$-rational points of $S$
is Zariski dense.

\ni b) Now let $C \subset S$ be a smooth curve of genus $1$, 
and assume that
$S$,
$\pi$ and
$C$ are defined over a field $K$. Then there is a finite 
extension $K'$ of $K$ such that the set $S(K')$ of
$K'$-rational points of
$S$ is Zariski dense.
\end{cor}

\

Part a) of Corollary~\ref{ellipticcor} follows from 
Theorem~\ref{ellipticlemma} as before: we see that for all but
finitely many points $p \in C(K)$, the fiber of 
$S \to \P^1$ over $\pi(p)$ has infinitely many rational points. As for
part b), we have to make an extension of our 
ground field $K$ simply to ensure that the curve $C$ has infinitely
many rational points, and then the argument 
proceeds as before. Of course we can drop the hypothesis in
Theorem~\ref{ellipticlemma} that $S$ is irrational.

\

\ni {\bf Proof of Theorem~\ref{ellipticlemma}}. 
The proof is analogous to that of Theorem~\ref{chowapproach}.
We begin by making a base change: we let $T$ be the incidence correspondence
$$ 
T \; = \; \bigl\{(p,q) : q \in E_{\pi(p)} \bigr\} \; \subset \; C \times S \, .
$$
As before,  $T$ is simply the fiber product
$$
T \; = \; C \times_{\P^1} S \, ;
$$
in other words, $T \to C$ is the fibration obtained 
by applying the base change $\pi|_C : C \to \P^1$  to the 
 fibration
$S
\to \P^1$. In particular,
$T$ is an $m$-sheeted cover of $S$, branched over 
the union of the fibers $E_\l$ of $S \to \P^1$ such
that $E_\l$ is tangent to $C$. Again, 
$T$ will have at worst isolated singularities, since by
the hypothesis that $S$ is smooth and all 
fibers of $\pi$ are irreducible it follows that all fibers are reduced as
well.

\

Note that $T \to C$ has a tautologous section
$$
\Sigma \; = \; \bigl\{(p,p) : p \in C \bigr\} \; \subset \; T \, .
$$
As a divisor, the pullback $\nu^{*}(C)$ of the curve $C$ under the
$m$-sheeted covering $\nu : T \to S$ is thus a sum
$$
\nu^{*}(L) \; = \; \Sigma + R
$$
with $R \subset T$ flat of degree $m-1$ over $C$. 
As before, since the  general fiber of $C$ over $\P^1$ is
reduced, $R$ does not contain $\Sigma$. 

\

Now let $\phi \in A_1(T)$ be the class of a fiber of $T \to C$, 
$\sigma$ the class of the section $\Sigma$ and
$\rho$ the class of $R$. The key ingredient in our proof is the

\

\begin{lm} The classes $\sigma$, $\rho$ and $\phi \in A_1(T)$ are
 independent in the group
$A_1(T)$ of Weil divisors mod linear equivalence on $T$.
\end{lm}

\

\ni {\bf Proof}. We calculate the matrix of intersection
products of the classes $\sigma$, $\rho$ and $\phi \in A_1(T)$. 
Three of these numbers are readily
calculated. To begin with, $\phi$ is the class of a fiber 
of the map $T \to L$, so of course $\phi^2=0$; and
since $\Sigma$ and $R$ meet each fiber in $1$ and $m-1$ 
points respectively, we have $(\phi \cdot \sigma) = 1$
and $(\phi \cdot \rho) = m-1$. 

As before, we don't know anything about $c = (\sigma \cdot \rho)$ except that it
is positive (if $T$ were smooth, if would have to be $2m-2$). 
Finally,  let $-b$ be
the self-intersection of the curve $C$ on $S$. By the 
hypothesis that $S$  is irrational, the canonical class $K_S$ is
a nonnegative (rational) multiple of the class of a 
fiber of $S \to \P^1$, and so has nonnegative intersection with
$C$; it follows that $b \geq 2$ if $C$ 
is rational, and $g \geq 0$ is $C$ has genus $1$.

\

Now, to calculate 
$\sigma^2$ and $\rho^2$, we use the relation
$$
\nu^*C \; = \; \Sigma + R \, .
$$
It follows that
$$
\sigma^2 \; = \; (\sigma \cdot [\nu^*C] - \rho) \, ,
$$
and since, by the push-pull formula,
$$
(\Sigma \cdot \nu^*C)_T \; = \; (\nu_*\Sigma \cdot C)_S  \; = \; (C \cdot C)_S  \; = \;  -b 
$$
we have
$$
\sigma^2 \; = \; -b -c \, .
$$
Similarly, from the same relation it follows that
$$
\rho^2 \; = \; (\rho \cdot [\nu^*C] - \sigma) \, .
$$
By the push-pull formula,
$$
(R \cdot \nu^*C)_T \; = \; (\nu_*R \cdot C)_S  \; = \; ((m-1)C \cdot C)_S  \; = \;  -(m-1)b 
$$
and so
$$
\rho^2 \; = \; -(m-1)b -c \, .
$$
In sum, then, we have the following table of intersection products
\begin{center}
\renewcommand{\arraystretch}{1.25}
\begin{tabular}{ l | c | c | c | }
& $\phi$ & $\sigma$ & $\rho$ \\
\hline
 $\phi$ & 0 & 1 & $m-1$ \\
\hline
$\sigma$ & 1 & $-b-c$ & $c$ \\
\hline
$\rho$ & $m-1$ & $c$ & $-(m-1)b-c$ \\
\hline
\end{tabular}
\end{center}

\

The determinant of this matrix is
$$
c(m-1) + (m-1)b + c + c(m-1) + (b+c)(m-1)^2
$$
and since $b$ is nonnegative and $c$ positive we are done.

\


\begin{thebibliography}{99}

\bibitem{barth-peters-vdv} W. Barth, C. Peters, A. van de Ven,
{\em Compact complex surfaces}, Ergebnisse der
Mathematik und ihrer Grenzgebiete, 
Springer-Verlag, Berlin-New York, (1984). 



\bibitem{beauville} A. Beauville, {\em    
Varietes de Prym et jacobiennes intermediaires},
Ann. Sci. Ecole Norm. Sup. (4) {\bf 10}, no. 3, (1977),  309--391.



\bibitem{beauville-3}
A. Beauville, {\em Varietes rationnelles et unirationnelles}, 
Algebraic geometry -- open problems (Ravello, 1982), LN in Math. {\bf 997}, 
Springer, Berlin-New York, (1983), 16--33.


\bibitem{clemens} C. H. Clemens, {\em Double solids}, 
Adv. Math. {\bf 47}, (1983), 107--230.

\bibitem{collino} A. Collino, {\em Lines on quartic threefolds}, 
J. London Math. 
Soc. (2) {\bf 19}, (1979), 257--267.

\bibitem{dem} V. A. Dem'janenko, {\em
Rational points of a class of algebraic curves}, 
A.M.S. Translations (2) {\bf 66},  (1968), 246--272.


\bibitem{FMT} J. Franke, Yu. I. Manin, Yu. Tschinkel, {\em
Rational points of bounded height on Fano varieties}, 
Invent. Math. {\bf 95}, no. 2, (1989), 421--435.

\bibitem{Fulton} W. Fulton, {\em
Intersection theory}, Ergebnisse der Mathematik und ihrer
Grenzgebiete (3), Springer-Verlag, Berlin-New
York, (1984).

\bibitem{isk}  V. A. Iskovskikh, {\em 
On the rationality problem for conic bundles}, Duke
Math. J. {\bf 54}, no. 2,  (1987), 271--294.

\bibitem{isk3} V.A.  Iskovskikh, {\em Fano threefolds. I},
Izv. Akad. Nauk SSSR Ser. Mat. {\bf 41}, no. 3, (1977), 516--562.

\bibitem{isk4} V.A.  Iskovskikh, 
{\em  Fano threefolds. II}, 
Izv. Akad. Nauk SSSR Ser. Mat. {\bf 42}, no. 3, (1978), 506--549. 

\bibitem{isk-manin} V.A.  Iskovskikh, Yu. I. Manin,
{\em Three dimensional quartics and counterexamples to L\"uroth problem},
Mat. Sbornik {\bf 86}, (1971), 140-166.

\bibitem{kollar-miyaoka-mori-92-2}
J. Koll{\'a}r, Y. Miyaoka, Sh. Mori,
{\em Rational curves on {F}ano varieties},
In {\em Classification of irregular varieties (Trento, 1990)},
LN in Math. {\bf 1515}, Springer, Berlin, (1992), 100--105. 

\bibitem{Manin} Yu. I. Manin, 
{\em Notes on the arithemetic of Fano threefolds},
Compositio Math. {\bf 85}, (1993), 37--55.

\bibitem{Manin2} Yu. I. Manin, 
{\em The $p$-torsion of elliptic curves is uniformly bounded},
Izv. Akad. Nauk. SSSR. {\bf 33}, (1969), 433--438.

\bibitem{Merel} L. Merel,  {\em 
Bornes pour la torsion des courbes elliptiques sur les corps de
nombres},  Invent. Math. {\bf 124}, no. 1-3, (1996), 437--449.

\bibitem{Mori-Mukai-1} Sh. Mori, Sh. Mukai, {\em  
On Fano $3$-folds with $B\sb{2}\geq 2$},  Algebraic
varieties and analytic varieties (Tokyo, 1981),
Adv. Stud. Pure Math. 1, North-Holland,
Amsterdam-New York, (1983), 101--129. 

\bibitem{Mori-Mukai-2} Sh. Mori, Sh. Mukai, {\em  
Classification of Fano $3$-folds with
$B_{2}\geq 2$},  Manuscripta Math. {\bf 36}, no. 2, (1981/82), 147--162. 


\bibitem{murre} J. P. Murre, {\em  
Classification of Fano threefolds according to Fano and Iskovskikh},
Algebraic threefolds (Varenna, 1981), LN in Math. {\bf  947}, 
Springer, Berlin-New York, (1982), 35--92.

\bibitem{noguchi} J. Noguchi, {\em A higher-dimensional analogue of Mordell's conjecture over function fields},
 Math. Ann. {\bf 258} (1981/82), no. 2, 207--212

\bibitem{shaf} I. I. Piatetski-Shapiro, I. Shafarevitch, {\em Torelli's theorem for algebraic surfaces of type K3},
(Russian) Izv. Akad. Nauk SSSR {\bf 35} (1971), 530--572

\bibitem{sil} J. H. Silverman, 
{\em Heights and the specialization map for families of abelian varieties},
J. reine und angew. Math. {\bf 342}, (1983), 197--211.


\bibitem{vojta} P. Vojta, {\em Diophantine approximation and value 
distribution theory}, LN in Math. {\bf 1239}, Springer, Berlin-New York, (1987).

\end{thebibliography}
\end{document}